\documentclass[12pt,twoside,leqno]{amsart}

\usepackage{amssymb,amsbsy,amsmath,amsfonts,eucal,amsthm,
amssymb,graphicx,times,xypic,color,wasysym}

\usepackage[T1]{fontenc} 
\definecolor{blue}{cmyk}{1.,1.,0.,0.53}
\definecolor{red}{cmyk}{0.,1.,1.,0.53}
\definecolor{green}{cmyk}{1.,0.,1.,0.53}

\newcommand{\C}{\mathbb{C}}\newcommand{\K}{\mathbb{K}}
\newcommand{\N}{\mathbb{N}}
\newcommand{\R}{\mathbb{R}}

\begin{document}

\title[Vanishing Hachtroudi-Chern CR curvature]{
Vanishing Hachtroudi curvature
\\
and local equivalence to the Heisenberg pseudosphere
}




\author{Jo\"el Merker}

\address{D\'epartement de math\'ematiques 
et applications, \'Ecole Normale Sup\'erieure, Paris}
\email{merker@dma.ens.fr} 

\date{\number\year-\number\month-\number\day. MSC:
32V40 (Primary). 35G20, 32W50, 58A15, 58A20 (Secondary)}

\begin{abstract}

To any completely integrable second-order system of 
real or complex partial differential equations:
\[
y_{x^{k_1}x^{k_2}}
=
F_{k_1,k_2}
\big(
x^1,\dots,x^n,\,y,\,y_{x^1},\dots,y_{x^n}
\big)
\]
with $1 \leqslant k_1,\, k_2 \leqslant n$ and with $F_{ k_1, k_2} =
F_{ k_2, k_1}$ in \underline{$n \geqslant 2$} independent variables
$(x^1, \dots, x^n)$ and in one dependent variable $y$, Mohsen
Hachtroudi associated in 1937 a normal projective (Cartan) connection,
and he computed its curvature. By means of a natural transfer of jet
polynomials to the associated submanifold of solutions, what the
vanishing of the Hachtroudi curvature gives can be precisely
translated in order to characterize when both families of Segre
varieties and of conjugate Segre varieties associated to a Levi
nondegenerate real analytic hypersurface $M$ in $\C^n$ ($n \geqslant
3$) can be straightened to be affine complex (conjugate) lines. In
continuation to a previous paper devoted to the quite distinct
$\C^2$-case, this then characterizes in an effective way those
hypersurfaces of $\C^{n+1}$ in higher complex dimension $n +
1\geqslant 3$ that are locally biholomorphic to a piece of the $(2n +
1)$-dimensional Heisenberg quadric, without any special assumption on
their defining equations.
\end{abstract}

\maketitle

\vspace{-1cm}

\begin{center}
\begin{minipage}[t]{11cm}
\baselineskip =0.35cm
{\scriptsize

\centerline{\bf Table of contents}

\smallskip

{\bf 1.~Introduction\dotfill 
\pageref{Section-1}.}

{\bf 2.~Segre varieties and differential equations\dotfill 
\pageref{Section-2}.}

{\bf 3.~Geometry of associated submanifolds of solutions\dotfill 
\pageref{Section-3}.}

{\bf 4.~Effective differential characterization of pseudosphericality
in $\C^{n+1}$\dotfill 
\pageref{Section-4}.}

}\end{minipage}
\end{center}

\section*{\S1.~Introduction}
\label{Section-1}

The explicit characterization of pseudosphericality of an arbitrary real
analytic local hypersurface sitting in the complex Euclidean space has
been (re)studied recently by Isaev in~\cite{ isa2009}, who employed
the famous Chern(-Moser) tensorial approach~\cite{ cm1974, ch1975} to
the concerned equivalence problem. But in the growing literature
devoted to Lie-group symmetries of Cauchy-Riemann manifolds, only a
very few articles underline that, already in his 1937
Ph.D. thesis~\cite{ ha1937} under the direction of his \'Elie
Cartan\,\,---\,\,who was around the same period also the master of
Chern\,\,---, the Iranian mathematician Mohsen Hachtroudi (cited
briefly only in~\cite{ ch1975}) constructed directly an {\em explicit}
normal projective Cartan connection canonically associated to any
completely integrable system of real or complex partial differential
equations:
\begin{equation}
\label{second-F}
y_{x^{k_1}x^{k_2}}
(x)
=
F_{k_1,k_2}
\big(
x^1,\dots,x^n,\,y,\,y_{x^1},\dots,y_{x^n}
\big)
\end{equation}
in \underline{$n \geqslant 2$} independent variables $x^1, \dots, x^n$
and in one dependent variable $y$, by endeavouring in a successful way
to generalize the celebrated paper~\cite{ ca1924}. Chern's clever
observation in 1974 that Hachtroudi's 37 years-old approach was
intrinsically related to the nascent higher-dimensional CR geometry
was followed, in his two papers in question, by his technical
contribution of redoing (only) {\em parts} of Hachtroudi's effective
computations, following the alternative (heavier, though essentially
equivalent) strategy of constructing {\em a posteriori} the projective
connection, after having reinterpreted at the beginning the problem in
terms of the wide and powerful {\sl Cartan Method of
Equivalence}. Thus, one should be aware, historically speaking, that
in the original reference~\cite{ ha1937}, much more complete geometric
and computational aspects were published long before, though they were
expressed in a purely analytic and somewhat elliptic language which,
unfortunately for us at present times, does not transmit in words and
with figures all the underlying geometric meanings which were clear
then to \'Elie Cartan.

Because Hachtroudi was able to write down explicitly his curvature
tensors, he deduced the second-order system~\thetag{
\ref{zero-hachtroudi}}\,\,---\,\,below\,\,---\,\,of partial
differential equations that the functions $F_{ k_1, k_2}$ should
satisfy in order that the system~\thetag{ \ref{second-F}} be
equivalent, through a point transformation $(x, y) \mapsto (x', y') =
\big( x' ( x, y), \, y' ( x, y) \big)$ to the simplest system: $y_{
{x'}^{ k_1} {x'}^{ k_2} }' ( x') = 0$, with all right-hand sides being
zero. In the present article, a companion and a follower of a
preceding one~\cite{ me2010c} devoted to the quite different
$\C^2$-case, we will apply, to the higher-dimensional characterization
of pseudosphericality, this effective necessary and sufficient
condition~\thetag{ \ref{zero-hachtroudi}} due to Hachtroudi which,
however and {\em inexplicably}, is totally inextant in the two
contributions of Chern.  We hope in this way to complete the explicit
characterization of pseudosphericality 
for rigid or even tube hypersurfaces
that was obtained recently by Isaev in~\cite{ isa2009}, because
apparently, the general (nonrigid) case was still open in the
specialized field.

\smallskip

We now start the exposition. Let $M$ be a local real analytic in $\C^{
n+1}$. Though the basic definitions, lemmas and propositions of the
theory are valid in any complex dimension $n+1 \geqslant 2$, there is
a strong computational difference between the two characterizations of
sphericality for $n = 1$ (compare~\cite{ me2010b}) and of
pseudosphericality for $n
\geqslant 2$ (presently), so that, in order to fix the ideas, it will
be assumed throughout the paper\,\,---\,\,and recalled when
necessary\,\,---\,\,that the CR dimension $n$ is always $\geqslant 2$.

Locally in a neighborhood of one of its points $p$, the hypersurface
$M$ may be represented, in any system of local holomorphic
coordinates:
\[
t=(w,z)
\in\C^n\times\C
\]
vanishing at $p$ for which the $w$-axis is not complex-tangent to $M$
at $p$, by a so-called {\em complex} defining
equation\,\,---\,\,Section~2 provides further
informations\,\,---\,\,of the form:
\begin{equation}
\label{Theta-defining}
w 
= 
\Theta\big(z,\,\overline{z},\overline{w})
=
\Theta\big(z,\,\overline{t}\big),
\end{equation}
or equivalently in a more expanded form which exhibits all the indices: 
\[
w
=
\Theta\big(z_1,\dots,z_n,\,
\overline{z}_1,\dots,\overline{z}_n,\,\overline{w}\big)
=
\Theta\big(z_1,\dots,z_n,\,
\overline{t}_1,\dots,\overline{t}_n,\overline{t}_{n+1}\big).
\]
Then $M$ localized at $p$ is called {\sl pseudospherical} (at $p$) if
it is biholomorphic to a piece of one Heisenberg pseudosphere:
\begin{equation}
\label{signature-levi-form}
{\rm Im}\,w'
=
\vert z_1'\vert^2
+\cdots+
\vert z_q'\vert^2
-
\vert z_{q+1}'\vert^2
-\cdots-
\vert z_n'\vert^2,
\end{equation}
for some $q$ with $0\leqslant q\leqslant n$, the number
of positive eigenvalues of the nondegenerate Levi form.
Next, let us introduce the following Jacobian-like determinant:  
\[
\Delta
:=
\left\vert
\begin{array}{cccc}
\Theta_{\overline{z}_1} & \cdots & \Theta_{\overline{z}_n} 
& \Theta_{\overline{w}}
\\
\Theta_{z_1\overline{z}_1} & \cdots & \Theta_{z_1\overline{z}_n} 
& \Theta_{z_1\overline{w}}
\\
\cdot\cdot & \cdots & \cdot\cdot & \cdot\cdot
\\
\Theta_{z_n\overline{z}_1} & \cdots & \Theta_{z_n\overline{z}_n} 
& \Theta_{z_n\overline{w}}
\end{array}
\right\vert
=
\left\vert
\begin{array}{cccc}
\Theta_{\overline{t}_1} & \cdots & \Theta_{\overline{t}_n} 
& \Theta_{\overline{t}_{n+1}}
\\
\Theta_{z_1\overline{t}_1} & \cdots & \Theta_{z_1\overline{t}_n} 
& \Theta_{z_1\overline{t}_{n+1}}
\\
\cdot\cdot & \cdots & \cdot\cdot & \cdot\cdot
\\
\Theta_{z_n\overline{t}_1} & \cdots & \Theta_{z_n\overline{t}_n} 
& \Theta_{z_n\overline{t}_{n+1}}
\end{array}
\right\vert.
\]
For any index $\mu \in \{ 1, \dots, n, n+1\}$ and for any index $\ell
\in \{ 1, \dots, n\}$, let also $\Delta_{ [ 0_{ 1+\ell}]}^\mu$ denote
the same determinant, but with its $\mu$-th column replaced by the
transpose of the line $(0 \cdots 1 \cdots 0)$ with $1$ at the
$(1+\ell)$-th place, and $0$ elsewhere, its other columns being
untouched. One easily convinces oneself (but see also Section~2) that
$M$ is Levi-nondegenerate at $p$\,\,---\,\,which is the origin of our
system of coordinates\,\,---\,\, if and only if $\Delta$ does not
vanish at the origin, whence $\Delta$
is nowhere zero in some sufficiently small
neighborhood of the origin.  Similarly, for any indices $\mu, \nu,
\tau \in \{ 1, \dots, n, n+1\}$, denote by $\Delta_{ [ \overline{
t}^\mu \overline{ t}^\nu]}^\tau$ the same determinant as $\Delta$, but
with only its $\tau$-th column replaced by the transpose of the line:
\[
\big(
\Theta_{\overline{t}^\mu\overline{t}^\nu}\ \ 
\Theta_{z_1\overline{t}^\mu\overline{t}^\nu}\ \
\cdots\ \
\Theta_{z_n\overline{t}^\mu\overline{t}^\nu}
\big),
\]
other columns being again untouched.  All these determinants $\Delta$,
$\Delta_{ [ 0_{ 1+\ell}]}^\mu$, $\Delta_{ [ \overline{ t}^\mu
\overline{ t}^\nu]}^\tau$ are visibly universal differential
expressions depending upon the second-order jet $J_{ z, \overline{ z},
\overline{ w}}^2 \Theta$ and upon the third-order jet $J_{ z, \overline{
z}, \overline{ w}}^3 \Theta$.

\smallskip\noindent{\bf Main Theorem.}  
{\em An arbitrary, not necessarily rigid, real analytic hypersurface
$M \subset \C^{ n+1}$ with \underline{$n \geqslant 2$} which is Levi
nondegenerate at one of its points $p$ and has a complex defining
equation of the form~\thetag{ \ref{Theta-defining}} in some system of
local holomorphic coordinates $t = (z, w) \in \C^n \times \C$
vanishing at $p$, is 
pseudospherical at $p$ {\em if and only if} its complex
graphing function $\Theta$ satisfies the following explicit nonlinear
fourth-order system of partial differential equations:}
\[
\boxed{
\scriptsize
\aligned
0
&
\equiv
\sum_{\mu=1}^{n+1}\,\sum_{\nu=1}^{n+1}
\bigg[
\Delta_{[0_{1+\ell_1}]}^\mu
\cdot
\Delta_{[0_{1+\ell_2}]}^\nu
\bigg\{
\Delta
\cdot
\frac{\partial^4\Theta}{
\partial z_{k_1}\partial z_{k_2}
\partial\overline{t}_\mu\partial\overline{t}_\nu}
-
\sum_{\tau=1}^{n+1}\,
\Delta_{[\overline{t}^\mu\overline{t}^\nu]}^\tau
\cdot
\frac{\partial^3\Theta}{
\partial z_{k_1}\partial z_{k_2}\partial\overline{t}^\tau}
\bigg\}
-
\\
&
-
{\textstyle{\frac{\delta_{k_1,\ell_1}}{n+2}}}\,
\sum_{\ell'=1}^n\,
\Delta_{[0_{1+\ell'}]}^\mu
\cdot
\Delta_{[0_{1+\ell_2}]}^\nu
\bigg\{
\Delta
\cdot
\frac{\partial^4\Theta}{
\partial z_{\ell'}\partial z_{k_2}
\partial\overline{t}_\mu\partial\overline{t}_\nu}
-
\sum_{\tau=1}^{n+1}\,
\Delta_{[\overline{t}^\mu\overline{t}^\nu]}^\tau
\cdot
\frac{\partial^3\Theta}{
\partial z_{\ell'}\partial z_{k_2}\partial\overline{t}^\tau}
\bigg\}
-
\\
&
-
{\textstyle{\frac{\delta_{k_1,\ell_2}}{n+2}}}\,
\sum_{\ell'=1}^n\,
\Delta_{[0_{1+\ell_1}]}^\mu
\cdot
\Delta_{[0_{1+\ell'}]}^\nu
\bigg\{
\Delta
\cdot
\frac{\partial^4\Theta}{
\partial z_{\ell'}\partial z_{k_2}
\partial\overline{t}_\mu\partial\overline{t}_\nu}
-
\sum_{\tau=1}^{n+1}\,
\Delta_{[\overline{t}^\mu\overline{t}^\nu]}^\tau
\cdot
\frac{\partial^3\Theta}{
\partial z_{\ell'}\partial z_{k_2}\partial\overline{t}^\tau}
\bigg\}
-
\\
&
-
{\textstyle{\frac{\delta_{k_2,\ell_1}}{n+2}}}\,
\sum_{\ell'=1}^n\,
\Delta_{[0_{1+\ell'}]}^\mu
\cdot
\Delta_{[0_{1+\ell_2}]}^\nu
\bigg\{
\Delta
\cdot
\frac{\partial^4\Theta}{
\partial z_{k_1}\partial z_{\ell'}
\partial\overline{t}_\mu\partial\overline{t}_\nu}
-
\sum_{\tau=1}^{n+1}\,
\Delta_{[\overline{t}^\mu\overline{t}^\nu]}^\tau
\cdot
\frac{\partial^3\Theta}{
\partial z_{k_1}\partial z_{\ell'}\partial\overline{t}^\tau}
\bigg\}
-
\\
&
-
{\textstyle{\frac{\delta_{k_2,\ell_2}}{n+2}}}\,
\sum_{\ell'=1}^n\,
\Delta_{[0_{1+\ell_1}]}^\mu
\cdot
\Delta_{[0_{1+\ell'}]}^\nu
\bigg\{
\Delta
\cdot
\frac{\partial^4\Theta}{
\partial z_{k_1}\partial z_{\ell'}
\partial\overline{t}_\mu\partial\overline{t}_\nu}
-
\sum_{\tau=1}^{n+1}\,
\Delta_{[\overline{t}^\mu\overline{t}^\nu]}^\tau
\cdot
\frac{\partial^3\Theta}{
\partial z_{k_1}\partial z_{\ell'}\partial\overline{t}^\tau}
\bigg\}
+
\\
&
\ \ \ \ \
+
{\textstyle{\frac{1}{(n+1)(n+2)}}}
\cdot
\big[
\delta_{k_1,\ell_1}\delta_{k_2,\ell_2}
+
\delta_{k_2,\ell_1}\delta_{k_1,\ell_2}
\big]
\cdot
\\
&
\ \ \ \ \
\cdot
\sum_{\ell'=1}^n\,\sum_{\ell''=1}^n\,
\Delta_{[0_{1+\ell'}]}^\mu
\cdot
\Delta_{[0_{1+\ell''}]}^\nu
\bigg\{
\Delta
\cdot
\frac{\partial^4\Theta}{
\partial z_{\ell'}\partial z_{\ell''}
\partial\overline{t}_\mu\partial\overline{t}_\nu}
-
\sum_{\tau=1}^{n+1}\,
\Delta_{[\overline{t}^\mu\overline{t}^\nu]}^\tau
\cdot
\frac{\partial^3\Theta}{
\partial z_{\ell'}\partial z_{\ell''}\partial\overline{t}^\tau}
\bigg\},
\endaligned}
\] 
{\em for all pairs of indices $(k_1, k_2)$ with $1 \leqslant k_1, k_2
\leqslant n$, and for all pairs of indices $(\ell_1, \ell_2)$ with $1
\leqslant \ell_1, \ell_2 \leqslant n$. 
}\medskip

The written system is effective: no implicit formal expression is
involved and pseudosphericality 
is characterized directly and only in terms
of $\Theta$.

\smallskip

Now, here is a summarized description of our arguments of proof.  A
bit similarly as for the $\C^2$-case\,\,---\,\,but with major
differences afterwards\,\,---\,\,which was already studied in~\cite{
me2010b}, we may associate to any such Levi nondegenerate real
analytic local hypersurface $M \subset \C^{ n+1}$ of equation $w =
\Theta ( z, \overline{ z}, \overline{ w})$ a uniquely defined
system of second-order partial differential equations:
\begin{equation}
\label{second-Phi}
w_{z_{k_1}z_{k_2}}(z)
=
\Phi_{k_1,k_2}
\big(
z,\,w(z),\,w_z(w)
\big)
\ \ \ \ \ \ \ \ \ \ \ \ \ \ \ \ \ \
{\scriptstyle{(1\,\leqslant\,k_1,\,\,k_2\,\leqslant\,n)}}
\end{equation}
with $\Phi_{ k_1, k_2} = \Phi_{ k_2, k_1}$, simply by eliminating the
two variables $\overline{ z}$ and $\overline{ w}$, viewed as
parameters, from the set of $n+1$ equations\footnote{\,
This process appears for instance in the references~\cite{ enlie1888,
ha1937, ch1975, su2001, su2002, bie2007, me2009}.
}: 
\[
w(z)=\Theta\big(z,\overline{z},\overline{w}\big),
\ \ \ \ \
w_{z_1}(z)=
{\textstyle{\frac{\partial\Theta}{\partial z_1}}}
\big(z,\overline{z},\overline{w}\big),
\ 
\dots\dots, \ \
w_{z_n}(z)=
{\textstyle{\frac{\partial\Theta}{\partial z_n}}}
\big(z,\overline{z},\overline{w}\big),
\]
---\,\,the assumption that the Jacobian determinant $\Delta$ is
nonvanishing at the origin being precisely the one which guarantees,
technically speaking, that the classical (holomorphic) implicit
function theorem applies\,\,---\,\,and then by replacing the so obtained
values for $\overline{ z}$ and $\overline{ w}$ in all second order
derivatives $\frac{ \partial^2 \Theta}{ \partial z_{ k_1} z_{ k_2}}
\big( z, \overline{ z}, \overline{ w}\big)$, {\em see}~\thetag{
  \ref{Phi-second}} below.  Trivially, this system is completely
integrable, for we just derived it from its general solution $w(z) :=
\Theta \big( z, \, \overline{ z}, \, \overline{ w} \big)$, where
$(\overline{ z}, \overline{ w})$ are understood as parameters.

As we said, Hachtroudi showed that the curvature of the projective
normal (Cartan) connection he associated with the system~\thetag{
\ref{second-F}} vanishes if and only if the right-hand side functions
$F_{ k_1, k_2}$ satisfy the following explicit differential system,
which is {\em linear} in terms of their second-order derivatives 
(all of which, notably, appear only with respect to the $y_{ x^\ell}$):
\begin{equation}
\label{zero-hachtroudi}
\footnotesize
\aligned
0
&
\equiv
\frac{\partial^2 F_{k_1,k_2}}{\partial y_{x^{\ell_1}}y_{x^{\ell_2}}}
-
\\
&
\ \ \ \ \
-
{\textstyle{\frac{1}{n+2}}}\,
\sum_{\ell'=1}^n
\left(
\delta_{k_1,\ell_1}
\frac{\partial^2F_{\ell',k_2}}{
\partial y_{x^{\ell'}}\partial y_{x^{\ell_2}}}
+
\delta_{k_1,\ell_2}
\frac{\partial^2F_{\ell',k_2}}{
\partial y_{x^{\ell_1}}\partial y_{x^{\ell'}}}
+
\delta_{k_2,\ell_1}
\frac{\partial^2F_{k_1,\ell'}}{
\partial y_{x^{\ell'}}\partial y_{x^{\ell_2}}}
+
\delta_{k_2,\ell_2}
\frac{\partial^2F_{k_1,\ell'}}{
\partial y_{x^{\ell_1}}\partial y_{x^{\ell'}}}
\right)
+
\\
&
\ \ \ \ \
+
{\textstyle{\frac{1}{
(n+1)(n+2)}}}
\big[
\delta_{k_1,\ell_1}\delta_{k_2,\ell_2}
+
\delta_{k_2,\ell_1}\delta_{k_1,\ell_2}
\big]
\sum_{\ell'=1}^n\,\sum_{\ell''=1}^n\,
\frac{\partial^2F_{\ell',\ell''}}{
\partial y_{x^{\ell'}}\partial y_{x^{\ell''}}}
\ \ \ \ \ \ \ \ \ \ \ \ \ \ \ \ \ \ \ \ \ \ \ \
\begin{array}{c}
{\scriptstyle{(1\,\leqslant\,k_1,\,\,k_2\,\leqslant\,n)}}
\\
{\scriptstyle{(1\,\leqslant\,\ell_1,\,\,\ell_2\,\leqslant\,n)}}
\end{array}.
\endaligned
\end{equation}
Hachtroudi also showed that this latter condition, better known
nowadays amongst the {\em Several Complex Variables} community as {\sl
vanishing of Chern(-Moser) curvature} to which it indeed amounts,
characterizes the local equivalence, through a point transformation
$(x, y) \mapsto (x', y') = \big( x' ( x, y), \, y' ( x, y) \big)$, to
the simplest system: $y_{ {x'}^{ k_1} {x'}^{ k_2} }' ( x') = 0$. We
then remind the semi-known fact that $M$ is 
pseudospherical if and only if
its associated second-order system~\thetag{ \ref{second-Phi}} is
equivalent, through a local biholomorphism $(z, w) \mapsto (z', w') =
\big( z' ( z, w), \, w' ( z, w) \big)$ fixing the origin, to the
simplest system $w_{ z_{ k_1}' z_{ k_2}'}' ( z') = 0$.  So we may
apply to the functions $\Phi_{ k_1, k_2}$ Hachtroudi's vanishing
curvature equations~\thetag{ \ref{zero-hachtroudi}}, but still, the
$\Phi_{ k_1, k_2}$ are not expressed in terms of $\Theta$, for they
were constructed by employing some unpleasant implicit functions when
solving above for $\overline{ z}$ and $\overline{ w}$. Fortunately,
here similarly as in~\cite{ me2010b}, we may apply the techniques of
computational differential algebra sketched in~\cite{ me2009} in order
to explicitly express any algebraic expressions in the second-order
jet of the $\Phi_{ k_1, k_2}$ in terms of the fourth-order jet of
$\Theta$, and the appropriate general equation which we shall need:
\[
\footnotesize
\aligned
\frac{\partial^2\Phi_{k_1,k_2}}{
\partial w_{z_{\ell_1}}\partial w_{z_{\ell_2}}}
&
=
\frac{1}{\Delta^3}\,
\sum_{\mu=1}^{n+1}\,\sum_{\nu=1}^{n+1}\,
\Delta_{[0_{1+\ell_1}]}^\mu
\!\cdot
\Delta_{[0_{1+\ell_2}]}^\nu
\left\{\,
\Delta\cdot
\frac{\partial^4\Theta}{
\partial z_{k_1}\partial z_{k_2}
\partial\overline{t}^\mu\partial\overline{t}^\nu}
-
\sum_{\tau=1}^{n+1}\,\,
\Delta_{[\overline{t}^\mu\overline{t}^\nu]}^\tau
\cdot
\frac{\partial^3\Theta}{
\partial z_{k_1}\partial z_{k_2}
\partial\overline{t}^\tau}
\right\}
\endaligned
\]
will be obtained in Section~4 below, after rather lengthy but
elementary calculations, parts of which are inspired from~\cite{
me2006}. It is now essentially clear how one obtains the (boxed) long
fourth-order differential equations stated in the theorem, but in any
case, some complete details will be provided at the very end of the
paper.

\smallskip

To conclude this extensive introduction which was designed for readers
wanting to quickly embrace the contents, we would like to draw the
attention on the work~\cite{ me2010c}, whose manual calculations where
finalized in manuscript form already in 2003\footnote{\,
At the conference {\scriptsize\sf\em Cauchy-Riemann Analysis and
Geometry} organized by Ingo Lieb and Gerd Schmalz at the {\sl
Max-Planck Institut} of Bonn, 22--27 September 2003, the author gave a
talk the title of which was ``{\em Explicit Chern-Moser tensors}''.
}, 
and which will soon confirm the above theorem by following another
route, {\em viz.} by calculating explicitly the so-called
Chern(-Moser) tensor differential forms,
which might interest some contemporary CR geometers
better than the (essentially equivalent) original
Cartan-Hachtroudi(-Tanaka) approach. 

\section*{\S2.~Segre varieties and differential equations}
\label{Section-2}

\subsection*{ Real analytic hypersurfaces in $\C^{n+1}$}
Let us therefore consider an arbitrary real analytic hypersurface $M$
in $\C^{ n+1}$ with $n \geqslant 2$, and let us localize it around one
of its points, say $p \in M$. Then there exist complex affine
coordinates:
\[
(z,w)
=
(z_1,\dots,z_n,w)
=
\big(x_1+
{\scriptstyle{\sqrt{-1}}}\,y_1,\dots,
x_n+{\scriptstyle{\sqrt{-1}}}\,y_n,\,u+iv\big)
=
\big(x+{\scriptstyle{\sqrt{-1}}}\,y,\,
u+{\scriptstyle{\sqrt{-1}}}\,v\big)
\]
vanishing at $p$ in which $T_p M = \{ u = 0\}$, so that $M$ is
represented in a neighborhood of $p$ by a graphed defining equation of
the form:
\[
u
=
\varphi(x,y,v)=
\varphi\big(x_1,\dots,x_n,y_1,\dots,y_n,v\big),
\]
where the real-valued function:
\[
\varphi 
= 
\varphi(x,y,v)
= 
\sum_{k\in\N^n,\,l\in\N^n,\,m\in\N\atop
\vert k\vert+\vert l\vert+m\geqslant 2}\,
\varphi_{k,l,m}\,x^ky^lv^m
\in
\R\big\{x,y,u\big\},
\]
which possesses entirely arbitrary real coefficients $\varphi_{
k,l,m}$, vanishes at the origin: $\varphi ( 0) = 0$, together with all
its first order derivatives: $0 = \partial_{ x^k} \varphi ( 0 ) =
\partial_{ y^l} \varphi ( 0) = \partial_v \varphi ( 0)$.  By simply
rewriting this initial real equation of $M$ as:
\[
{\textstyle{\frac{w+\overline{w}}{2}}}
=
\varphi
\big(
{\textstyle{\frac{z+\overline{z}}{2}}},\,
{\textstyle{\frac{z-\overline{z}}{2{\scriptstyle{\sqrt{-1}}}}}},\,
{\textstyle{\frac{w-\overline{w}}{2{\scriptstyle{\sqrt{-1}}}}}}
\big),
\]
and then by solving the so written equation with respect to $w$, one
obtains an equation of the shape:
\[
w
=
\Theta\big(z,\,\overline{z},\,\overline{w}\big)
=
\sum_{k\in\N^n,\,l\in\N^n,\,m\,\in\,\N\atop
\vert k\vert+\vert l\vert+m\geqslant 1}\,
\Theta_{k,l,m}\,
z^k\,\overline{z}^l\,\overline{w}^m
\in
\C\big\{\overline{z},\,z,\,w\big\},
\]
whose right-hand side converges of course near the origin $(0, 0,
0) \in \C^n \times \C^n \times \C$ and whose coefficients $\Theta_{ k,
l, m} \in \C$ are {\em complex}.  Since $d\varphi ( 0) = 0$, one has
$\Theta = - \overline{ w} + {\sf order}\,2\,{\sf terms}$.

The paradox that any such {\em complex} equation provides in fact {\em
two} real defining equations for the {\em real} hypersurface $M$ which
is {\em one}-codimensional, and also in addition the fact that one
could as well have chosen to solve the above equation with respect to
$\overline{ w}$, instead of $w$, these two apparent ``contradictions''
are corrected by means of a fundamental, elementary statement that
transfers to $\Theta$ (in a natural way) the condition of reality:
\[
\overline{\varphi(x,y,u)}
=
\sum_{\vert k\vert+\vert l\vert+m\geqslant 1}\,
\overline{\varphi_{k,l,m}}\,
\overline{x}^k\overline{y}^l\overline{v}^m
=
\sum_{\vert k\vert+\vert l\vert+m\geqslant 1}\,
\varphi_{k,l,m}\,x^ky^lv^m
=
\varphi(x,y,v)
\]
enjoyed by the initial definining function $\varphi$.  In the sequel,
we shall work exclusively with $\Theta$;
the reader is referred to~\cite{ me2010b} for
justifications and 
motivations. 

\smallskip\noindent{\bf Theorem.}
(\cite{ mp2006}, p.~19)
{\em The complex analytic function $\Theta = \Theta ( z, \overline{
z}, \overline{ w} )$ with $\Theta = - \overline{ w} + {\sf O}(2)$
together with its complex conjugate:
\[
\overline{\Theta}
=
\overline{\Theta}\big(\overline{z},z,w)
=
\sum_{k\in\N^n,\,l\in\N^n,\,m\in\N}\,
\overline{\Theta}_{k,l,m}\,
\overline{z}^k\,z^l\,w^m
\in
\C\big\{\overline{z},\,z,\,w\big\}
\]
satisfy the two {\em (}equivalent by conjugation{\em )} functional
equations{\em :}}
\begin{equation}
\label{reality-Theta}
\aligned
\overline{w}\equiv 
& \
\overline{\Theta}
\big(\overline{z},z,\Theta(z,\overline{z},\overline{w})\big), 
\\
w\equiv 
& \
\Theta\big(z,\overline{z}, 
\overline{\Theta}(\overline{z},z,w)\big).
\endaligned
\end{equation}
{\em 
Conversely, given a local holomorphic function $\Theta ( z, \overline{
z}, \overline{ w} ) \in \C \{ z, \overline{ z}, \overline{ w} \}$,
$\Theta = - \, \overline{ w} + {\sf O} ( 2)$ which, in conjunction
with its conjugate $\overline{ \Theta} ( \overline{ z}, z, w)$,
satisfies this pair of equivalent identities, then the two zero-sets:}
\[
\big\{
0=-\,w+\Theta\big(z,\,\overline{z},\,\overline{w}\big)
\big\}
\ \ \ \ \ \ \ \ \ \ 
\text{\em and}
\ \ \ \ \ \ \ \ \ \ 
\big\{
0=-\,\overline{w}+
\overline{\Theta}\big(\overline{z},\,z,\,w\big)
\big\}
\] 
{\em coincide and define a local {\em one-codimensional} real analytic
hypersurface $M$ passing through the origin in $\C^{ n+1}$.}\medskip

\subsection*{ Levi nondegeneracy} Within the hierarchy 
of nondegeneracy conditions for real hypersurfaces initiated by
Diederich and Webster (\cite{ dw1980}, {\em see} also~\cite{ me2005a,
me2005b} for generalizations and a unification), Levi nondegeneracy is
the most studied. The classical definition may be found in~\cite{
bo1991} and in the survey of Chirka~\cite{ ch1991}, but the following
basic equivalent characterization can also be understood as a
definition in the present paper. One may show (\cite{ me2005a, 
me2005b, mp2006}) that it is biholomorphically invariant. 

\smallskip\noindent{\bf Lemma.}
(\cite{ mp2006}, p.~28)
{\em 
The real analytic hypersurface $M \subset \C^{ n+1}$ with $0 \in M$
represented in coordinates $(z_1, \dots, z_n, w)$ by a complex
defining equation of the form $w = \Theta ( z, \overline{ z},
\overline{ w})$ is Levi nondegenerate at the origin if and only if the
map:
\[
\big(\overline{z}_1,\dots,\overline{z}_n,\overline{w}\big)
\longmapsto
\left(
\Theta\big(0,\overline{z},\overline{w}\big),\,\,
{\textstyle{\frac{\partial\Theta}{\partial z_1}}}
\big(0,\overline{z},\overline{w}\big),
\dots,
{\textstyle{\frac{\partial\Theta}{\partial z_n}}}
\big(0,\overline{z},\overline{w}\big)
\right)
\]
has nonvanishing $(n+1) \times (n+1)$ Jacobian determinant at
$(\overline{ z}, \overline{ w}) = (0, 0)$.
}\medskip

It follows then that this Jacobian determinant, not restricted
to the origin: 
\begin{equation}
\label{Delta-nonzero}
\Delta
=
\Delta\big(z,\overline{z},\overline{w}\big)
:=
\left\vert
\begin{array}{cccc}
\Theta_{\overline{z}_1} & \cdots & \Theta_{\overline{z}_n} 
& \Theta_{\overline{w}}
\\
\Theta_{z_1\overline{z}_1} & \cdots & \Theta_{z_1\overline{z}_n} 
& \Theta_{z_1\overline{w}}
\\
\cdot\cdot & \cdots & \cdot\cdot & \cdot\cdot
\\
\Theta_{z_n\overline{z}_1} & \cdots & \Theta_{z_n\overline{z}_n} 
& \Theta_{z_n\overline{w}}
\end{array}
\right\vert
\end{equation}
does not vanish in some small neighborhood of the origin in $\C^n
\times \C^n \times \C$. Levi nondegeneracy at the central point, {\em
i.e.}  $\Delta \neq 0$ locally, will be assumed throughout
the present paper. 

\subsection*{ Associated system of partial differential equations}
At least since the publication in 1888 by Lie and Engel in Leipzig of
the {\em Theorie der Transformationsguppen}, it is known in a very
general context\,\,---\,\,{\em see} Chapter~10 of~\cite{ enlie1888}
and also~\cite{ seg1931, ha1937, ch1975, fa1980, me2009, bie2007,
me2010b}\,\,---\,\,that, to the whole family of {\sl Segre varieties}:
\[
S_{\overline{z},\overline{w}}
:=
\big\{
(z,w)\in\C^n\times\C\colon\,\,
w=\Theta\big(z,\overline{z},\overline{w}\big)
\big\}
\]
parametrized by the $n+1$ antiholomorphic variables $\big( \overline{
z}_1, \dots, \overline{ z}_n, \overline{ w} \big)$, one may
canonically associate a completely integrable second-order
system of partial differential equations whose general solution is
precisely the function $\Theta \big( z, \overline{ z}, \overline{ w}
\big)$. Indeed, considering $w$ as a function $w = w(z)$ of $(z_1,
\dots, z_n)$ in the defining equation of $M$, one differentiates it
once with respect to each variable $z_1, \dots, z_n$ so that one gets
the $n+1$ equations:
\[
w(z)=\Theta\big(z,\overline{z},\overline{w}\big),
\ \ \ \ \
w_{z_1}(z)=
{\textstyle{\frac{\partial\Theta}{\partial z_1}}}
\big(z,\overline{z},\overline{w}\big),
\ 
\dots\dots, \ \
w_{z_n}(z)=
{\textstyle{\frac{\partial\Theta}{\partial z_n}}}
\big(z,\overline{z},\overline{w}\big).
\]
Then by means of the implicit function theorem\,\,---\,\,which
applies precisely
thanks to the nonvanishing of $\Delta$\,\,---, one may clearly solve
for the $n+1$ antiholomorphic ``parameters'' $(\overline{ z},
\overline{ w})$, and this procedure provides a representation:
\[
\overline{z}_1
=
\zeta_1\big(z,\,w(z),\,w_z(z)\big),
\
\dots,
\ \
\overline{z}_n
=
\zeta_n\big(z,\,w(z),\,w_z(z)\big),
\ \ 
\overline{w}
=
\xi\big(z,\,w(z),\,w_z(z)\big)
\]
with certain $n+1$ uniquely defined local complex analytic functions
$\zeta_i ( z, w, w_z)$ and $\xi ( z, w, w_z)$ of $2n+1$ complex
variables. Utilizing these functions, one is then pushed to replace
$\overline{ z}$ and $\overline{ w}$ in all possible second-order
derivative:
\begin{equation}
\label{Phi-second}
\aligned
w_{z_{k_1}z_{k_2}}(z)
&
=
{\textstyle{\frac{\partial^2\Theta}{\partial z_{k_1}\partial z_{k_2}}}}
\big(z,\,\overline{z},\,\overline{w}\big)
\\
&
=
{\textstyle{\frac{\partial^2\Theta}{\partial z_{k_1}\partial z_{k_2}}}}
\big(z,\,
\zeta\big(z,\,w(z),\,w_z(z)\big),\,\,
\xi\big(z,\,w(z),\,w_z(z)\big)\big)
\\
&
=:
\Phi_{k_1,k_2}\big(z,\,w(z),\,w_z(z)\big)
\ \ \ \ \ \ \ \ \ \ \ \ \
{\scriptstyle{(k_1,\,k_2\,=\,1\,\cdots\,n)}},
\endaligned
\end{equation}
and this defines without ambiguity the associated system of partial
differential equations. It is of second order. It is {\sl complete}:
all second-order derivatives are functions of derivatives of lower
order $\leqslant 1$.  In a sense to be precised right now, it is also
{\sl completely integrable} because by construction, its general
solution is $\Theta \big( z, \overline{ z}, \overline{ w}\big)$.

\subsection*{ Geometric characterization of pseudosphericality}
It is well known that the unit {\em sphere}:
\[
S^{2n+1}
=
\big\{
(z_1,\dots,z_n,w)\in\C^n\times\C\colon\,\,
\vert z_1\vert^2+\cdots+\vert z_n\vert^2+\vert w\vert^2
=
1
\big\}
\]
in $\C^n$ minus one of its points, for instance: $S^{2n+1} \setminus
\{ p_\infty \}$ with $p_\infty := (0, \dots, 0, -1)$, is
biholomorphic, through the so-called {\sl Cayley transform}:
\[
(z_1,\dots,z_n,w)\longmapsto
\big(
{\textstyle{\frac{i\,z_1}{1+w}}},\,\,
\dots,\,\,
{\textstyle{\frac{i\,z_n}{1+w}}},\,\,
{\textstyle{\frac{1-w}{2+2\,w}}}
\big)
=:
(z_1',\dots,z_n',w')
\]
having inverse:
\[
(z_1',\dots,z_n',w')\longmapsto
\big(
{\textstyle{\frac{-2iz_1'}{1+2w'}}},\,\,
\dots,
{\textstyle{\frac{-2iz_n'}{1+2w'}}},\,\,
{\textstyle{\frac{1-2w'}{1+2\,w'}}}
\big)
=
(z_1,\dots,z_n,w)
\]
to the so-called {\sl standard Heisenberg sphere} of equation:
\[
w'
=
-\,\overline{w}'
+
z_1'\overline{z}_1'
+\cdots+
z_n'\overline{z}_n'
\]
which sits in the target space $(z', w')$. Hence in the particular
case when the Levi form of $M$ has only positive eigenvalues, namely
when $q = n$ in~\thetag{ \ref{signature-levi-form}}, it follows
clearly that $M$ is {\sl spherical} in the sense given in the
Introduction if and only if there exists a nonempty open neighborhood
$U_0$ of $0$ in $\C^{ n+1}$ such that $M \cap U_0$ is biholomorphic to
a piece of the unit sphere.  In general, there are $n-q$ negative
eigenvalues in the Levi form, and this justifies 
adding a ``{\em pseudo}''.

\smallskip\noindent{\bf Proposition.}
\label{characterization-sphericality}
{\em A Levi nondegenerate local real analytic hypersurface $M$ in
$\C^{ n+1}$ is locally biholomorphic to a piece of the Heisenberg
pseudosphere (hence pseudospherical) if and only if its associated
second-order ordinary complex differential equation is locally
equivalent to the second-order system:}
\[
w_{z_{k_1}'z_{k_2}'}'(z')
=
0
\ \ \ \ \ \ \ \ \ \ \ \ \
{\scriptstyle{(1\,\leqslant\,k_1,\,\,k_2\,\leqslant\,n)}},
\]
{\em with identically vanishing right-hand side.}

\proof
The $n=1$ case, treated in great details by a previous
reference~\cite{ me2010b}, generalizes here with rather evident
adaptations, hence will be skipped. As $n \geqslant 2$ throughout the
present paper, one may also argue by slicing $\C^{ n+1}$ by all
possible copies of $\C^2$ which pass through the origin and which
contain the $w$-axis, so as to be able to apply the alreaday detailed
$n = 1$ case.
\endproof

Geometrically, the local equivalence of $M$ to the Heisenberg 
pseudosphere
means that, through some suitable local biholomorphism $(z, w) \mapsto
(z', w')$ fixing the origin, both its Segre varieties and its
conjugate Segre varieties (\cite{ me2005a, me2005b, mp2006}):
\[
S_{\overline{z},\overline{w}}
:=
\big\{
(z,w)\colon
w=\Theta\big(z,\overline{z},\overline{w}\big)
\big\}
\ \ \ \ \ \ \ \ \ \ 
\text{\rm and}
\ \ \ \ \ \ \ \ \ \ 
\overline{S}_{z,w}
:=
\big\{
(\overline{z},\overline{w})\colon
\overline{w}
=
\overline{\Theta}\big(\overline{z},z,w\big)
\big\}
\]
are mapped to the Segre and conjugate Segre
varieties of the Heisenberg pseudosphere: 
\[
S_{\overline{z}',\overline{w}'}'
=
\big\{
w'
=
-\,\overline{w}'
+
z'\overline{z}'
\big\}
\ \ \ \ \ \ \ \ \ \ 
\text{\rm and}
\ \ \ \ \ \ \ \ \ \ 
\big\{
\overline{w}'
=
-\,w'
+
\overline{z}'z'
\big\}
\]
which, visibly, are plain complex affine lines. 

\section*{\S3.~Geometry of associated submanifolds of solutions}
\label{Section-3}

\subsection*{ Completely integrable systems of
partial differential equations} 
The characterization of pseudosphericality we are dealing with holds in a
context more general than just CR geometry\footnote{\,
We will be very brief here, the reader being referred to~\cite{
me2009, me2010b} for the general theoretical considerations.
}. 
Accordingly, let $\K$
denote either the field $\C$ of complex numbers or the field $\R$ of
real numbers, let $x = (x^1, \dots, x^n) \in \K^n$ with again $n
\geqslant 2$\,\,---\,\,since the case $n = 1$ was already studied
in~\cite{ me2010b}\,\,---, let $y \in \K$, and consider a system of
the form~\thetag{ \ref{second-F}}. We will assume that it is {\sl
completely integrable} in the sense that the natural commutativity of
partial derivatives enjoyed trivially by the left-hand sides:
\[
\aligned
&
\partial^2y_{x^{k_1}x^{k_2}} 
\big/ 
\partial y_{ x^{k_3}}
=
\partial^2y_{x^{k_1}x^{k_3}} 
\big/ 
\partial y_{ x^{k_2}}
\\
&
\ \ \ \ \ \ \ \ \ \ \ \ \ \ \ \ \
{\scriptstyle{(1\,\leqslant\,k_1,\,\,k_2,\,\,k_3\,\leqslant\,n)}}
\endaligned
\]
imposes immediately to the right-hand side functions $F_{ k_1, k_2}$
that they satisfy the so-called {\sl compatibility conditions}:
\[
{\sf D}_{k_3}
\big(F_{k_1,k_2}\big)
=
{\sf D}_{k_2}
\big(F_{k_1,k_3}\big),
\]
where we have introduced the following $n$ {\sl total
differentiation} operators: 
\[
\aligned
{\sf D}_k
:=
{\textstyle{\frac{\partial}{\partial x^k}}}
&
+
y_{x^k}\,
{\textstyle{\frac{\partial}{\partial y}}}
+
\sum_{\ell=1}^n\,
F_{k,l}\,
{\textstyle{\frac{\partial}{\partial y_{x^\ell}}}}
\\
&
\ \ \ \ \ 
{\scriptstyle{(1\,\leqslant\,k\,\leqslant\,n)}}
\endaligned
\]
living on the first-order jet space $(x^1, \dots, x^n, y, y_{x^1},
\dots, y_{ x^n})$. One verifies that these compatibility conditions
amount to the fact that the $n$-dimensional tangential distribution
spanned by ${\sf D}_1, \dots, {\sf D}_n$ in the $(2n+1)$-dimensional
first-order jet space satisfies the classical Frobenius integrability
condition $\big[ {\sf D}_{ k'}, \, {\sf D}_{ k''} \big] = 0$, and then the
Clebsch-Frobenius theorem tells us that this distribution comes from a
local foliation by $n$-dimensional manifolds graphed over the
$x$-space that are naturally parametrized by $n+1$ auxiliary constants
(transversal directions)\,\,---\,\,call 
them $a^1, \dots, a^n, b \in \K$\,\,---,
namely the leaves of this local foliation may be explicitly
represented as sets of the shape:
\[
\aligned
\Big\{
\big(
x^1,\dots,x^n,\,
&
Q\big(x^1,\dots,x^n,\,a^1,\dots,a^n,b\big),
\\
&
S^1\big(x^1,\dots,x^n,\,a^1,\dots,a^n,b\big),\dots,
S^n\big(x^1,\dots,x^n,\,a^1,\dots,a^n,b\big)
\big)
\Big\},
\endaligned
\]
where $x^1, \dots, x^n$ vary freely and where $Q$, $S^1, \dots, S^n$
are certain graphing functions. In fact, the functions
$S^k$ are the first-order derivatives:
\[
S^1
=
Q_{x^1},\,\,
\dots\dots,\,\,
S^n
=
Q_{x^n}
\]
of the function $Q$, because by definition the integral curves of
every vector field ${\sf D}_k$ must be contained in such leaves, so
that one has:
\[
{\textstyle{\frac{\partial Q}{\partial x^k}}}
=
y_{x^k}\big\vert_{\text{\scriptsize\sf any leaf}}
=
S^k
\]
and furthermore also:
\[
{\textstyle{\frac{\partial S^l}{\partial x^k}}}
=
F^{k,l}\big\vert_{\text{\scriptsize\sf any leaf}},
\]
whence we see that the {\sl fundamental graphing function} $Q = Q ( x,
a, b)$ happens to be the {\em general solution} to the initially given
system of partial differential equations:
\[
\aligned
Q_{x^{k_1}x^{k_2}}(x,a,b)
&
\equiv
F_{k_1,k_2}
\big(
x,\,Q(x,a,b),\,Q_{x^1}(x,a,b),\dots,Q_{x^n}(x,a,b)
\big)
\\
&
\ \ \ \ \ \ \ \ \ \ \ \ \ \ \ \ \ \
{\scriptstyle{(1\,\leqslant\,k_1,\,\,k_2\,\leqslant\,n)}}.
\endaligned
\]
In the CR case, the fundamental function which is the general solution
to the associated system of partial diffential equations~\thetag{
\ref{second-Phi}} is obviously the complex defining function $\Theta
\big( z, \overline{ z}, \overline{ w} \big)$, where the $n+1$
quantities $(\overline{ z}, \overline{ w})$, viewed as independent
variables, play the role of the constants $(a, b)$.

As in the $n = 1$ case, the constants $(a^1, \dots, a^n, b)$ are best
interpreted as a set of $n+1$ {\sl initial conditions} $\big( y_{ x^1}
( 0), \dots, y_{ x^n} ( 0), - y(0) \big)$ or {\sl integration
constants}, so that we can assume without loss of generality that the
first-order terms in the fundamental function $Q$ are\footnote{\,
We put a minus sign in front of $y ( 0)$ so as to match up with our
choice of complex defining equation $w = - \, \overline{ w} + {\sf O}
( 2)$.
}: 
\[
Q(x,a,b)
=
-\,b
+
x^1a^1+\cdots+x^na^n
+
{\sf O}(\vert x\vert^2). 
\]
It is then clear that the map: 
\begin{equation}
\label{x-jet}
\aligned
\big(a^1,\dots,a^n,b\big)
&
\longmapsto
\big(
Q(0,a,b),\,
Q_{x^1}(0,a,b),\dots,
Q_{x^n}(0,a,b)
\big)
\\
&
\ \ \
=
\big(-b,a^1,\dots,a^n\big)
\endaligned
\end{equation}
is of rank $n+1$ at the origin, and this property remains also true
whatever one chooses as a fundamental function $Q ( x, a, b)$, that is
to say, without necessarily assuming it to be normalized as above,
which amounts to saying that\footnote{\,
Much more theoretical information is provided in~\cite{ me2009}.
},
in the parameter $(a,b)$-space, everything holds invariantly up to any
local $\K$-analytic transformation $(a, b) \mapsto (a', b')$ which
does {\em not} involve the variables $(x, y)$.

The way how one recovers the system of partial differential equations
is very similar to what we did in the CR case~\thetag{
\ref{second-Phi}}. Suppose indeed a bit more generally that we are
given any local $\K$-analytic function $Q = Q ( x, a, b)$ having the
property that its first-order $x$-jet map~\thetag{ \ref{x-jet}} is of
rank $n+1$ at $(a, b) = (0, 0)$.  Then in the $n+1$ equations:
\[
y(x)
=
Q(x,a,b),
\ \ \ 
y_{x^1}
=
Q_{x^1}(x,a,b),\,\dots\dots,\,\,
Q_{x^n}(x,a,b),
\]
we can solve, by means of the implicit function theorem, for the $n+1$
constants $(a^1, \dots, a^n, b)$, and this yields a representation:
\[
\aligned
a^k
&
=
A^k\big(x^1,\dots,x^n,y,y_{x^1},\dots,y_{x^n}\big)
\\
b
&
=
B\big(x^1,\dots,x^n,y,y_{x^1},\dots,y_{x^n}\big)
\endaligned
\]
for certain functions $A^1, \dots, A^n, B$ of $(2n+1)$ variables. Then
by replacing these obtained values for the $a^k$ and for $b$ in all
the possible second-order derivatives:
\[
\aligned
y_{x^{k_1}x^{k_2}}
&
=
Q_{x^{k_1}x^{k_2}}(x,a,b)
\\
&
=
Q_{x^{k_1}x^{k_2}}
\big(x,A(x,y,y_x),B(x,y,y_x)\big)
\\
&
=
F_{k_1,k_2}\big(x,y,y_x\big)
\endaligned
\]
it is rigorously clear that one may only recover the functions $F_{
k_1, k_2}$ we started with.

\section*{\S4.~Effective differential characterization 
\\
of pseudosphericality in $\C^{n+1}$}
\label{Section-4}

The $2n+1$ coordinates of the transformation considered at
the moment: 
\begin{equation}
\label{x-y-a-b}
\big(x^1,\dots,x^n,y,y_{x^1},\dots,y_{x^n}\big)
\longmapsto
\big(x^1,\dots,x^n,a^1,\dots,a^n,b\big)
\end{equation}
and those of its inverse are given by the collection of functions:
\[
\small
\left[
\aligned
x^j
&
=
x^j
\\
a^k
&
=
A^k\big(x^1,\dots,x^n,y,y_{x^1},\dots,y_{x^n}\big)
\\
b
&
=
B\big(x^1,\dots,x^n,y,y_{x^1},\dots,y_{x^n}\big)
\endaligned\right.
\ \ \ \ \
\text{\rm and}
\ \ \ \ \
\left[
\aligned
x^j
&
=
x^j
\\
y
&
=
Q\big(x^1,\dots,x^n,a^1,\dots,a^n,b\big)
\\
y_{x^k}
&
=
Q_{x^k}(x^1,\dots,x^n,a^1,\dots,a^n,b\big).
\endaligned\right.
\]
For uniformity and harmony, we shall admit by convention the
equivalences of notation:
\[
b
\equiv
a^{n+1}
\ \ \ \ \ \ \ \ \ \ 
\text{\rm and}
\ \ \ \ \ \ \ \ \ \
B
\equiv
A^{n+1}.
\]
Then by differentiating with respect to $y_{ x^\ell}$ each one of the
following $n+1$ identically satisfied equations:
\[
\small
\aligned
y
&
\equiv
Q
\big(
x^1,
\dots,x^n,\,
A^1\big(x^1,\dots,x^n,y,y_{x^1},\dots,y_{x^n}\big),
\dots,
\\
&
\ \ \ \ \ \ \ \ \ \ \
A^n\big(x^1,\dots,x^n,y,y_{x^1},\dots,y_{x^n}\big),
A^{n+1}\big(x^1,\dots,x^n,y,y_{x^1},\dots,y_{x^n}\big)
\big)
\\
y_{x^k}
&
\equiv
Q_{x^k}
\big(
x^1,
\dots,x^n,\,
A^1\big(x^1,\dots,x^n,y,y_{x^1},\dots,y_{x^n}\big),
\dots,
\\
&
\ \ \ \ \ \ \ \ \ \ \ \ \ \
A^n\big(x^1,\dots,x^n,y,y_{x^1},\dots,y_{x^n}\big),
A^{n+1}\big(x^1,\dots,x^n,y,y_{x^1},\dots,y_{x^n}\big)
\big),
\endaligned
\]
we get the following $n + n^2$ equations: 
\[
\aligned
0
&
\equiv
Q_{a^1}
{\textstyle{\frac{\partial A^1}{\partial y_{x^\ell}}}}
+\cdots+
Q_{a^n}
{\textstyle{\frac{\partial A^n}{\partial y_{x^\ell}}}}
+
Q_{a^{n+1}}
{\textstyle{\frac{\partial A^{n+1}}{\partial y_{x^\ell}}}}
\\
\delta_{k,\ell}
&
=
Q_{x^k a^1}
{\textstyle{\frac{\partial A^1}{\partial y_{x^\ell}}}}
+\cdots+
Q_{x^k a^n}
{\textstyle{\frac{\partial A^n}{\partial y_{x^\ell}}}}
+
Q_{x^k a^{n+1}}
{\textstyle{\frac{\partial A^{n+1}}{\partial y_{x^\ell}}}}
\\
&
\ \ \ \ \ \ \ \ \ \ \ \ \ \ \ \ \ \ \ \ \ \ \ \ \ \
{\scriptstyle{(k,\,\ell\,=\,1\,\cdots\,n)}}.
\endaligned
\]
Fixing any $\ell \in \{ 1, \dots, n\}$, thanks to the 
assumption (Levi nondegeneracy) that the
Jacobian determinant:
\[
\square
=
\square\big(a^1\vert\cdots\vert a^n\vert a^{n+1}\big)
:=
\left\vert
\begin{array}{cccc}
Q_{a^1} & \cdots & Q_{a^n} & Q_{a^{n+1}}
\\
Q_{x^1a^1} & \cdots & Q_{x^1a^n} & Q_{x^1a^{n+1}}
\\
\vdots & \ddots & \vdots & \vdots
\\
Q_{x^na^1} & \cdots & Q_{x^na^n} & Q_{x^na^{n+1}}
\end{array}
\right\vert,
\]
does not vanish, we may solve\,\,---\,\,just by means of Cramer's
rule\,\,---\,\,for the $n+1$ unknowns ${\frac{ \partial A^\mu}{
\partial y_{ x^\ell}}}$, the above system of $n+1$ equations, and this
gives us:
\begin{equation}
\label{A-square}
\frac{\partial A^\mu}{\partial y_{x^\ell}}
=
\frac{\square_{[0_{1+\ell}^\mu]}}{\square}
:=
\frac{
\square(a_1\vert\cdots\vert a^{\mu-1}\vert 0_{1+\ell}\vert
a^{\mu+1}\vert\cdots\vert a^{n+1})}{
\square(a^1\vert\cdots\vert a^{\mu-1}\vert a^\mu\vert a^{\mu+1}\vert
\cdots\vert a^{n+1})}, 
\end{equation}
where $0_{ 1+\ell}$ is a specific notation to denote the column
consisting of $n+1$ zeros piled up, except at the $(1+\ell)$-th level
from its top, where instead of $0$, one reads $1$, and where, as our
notation with vertical bars helps to guess:
\[
\small
\aligned
\square_{[0_{1+\ell}]}^\mu
=
\square(a_1\vert\cdots\vert a^{\mu-1}\vert^\mu\,0_{1+\ell}\vert
a^{\mu+1}\vert\cdots\vert a^{n+1})
:=
\endaligned
\]
\[
\small
\aligned
:=
\left\vert
\begin{array}{ccccccc}
Q_{a^1} & \cdots & Q_{a^{\mu-1}} & 0 & Q_{a^{\mu+1}} 
& \cdots & Q_{a^{n+1}}
\\
Q_{x^1a^1} & \cdots & Q_{x^1a^{\mu-1}} & 0 & Q_{x^1a^{\mu+1}}
& \cdots & Q_{x^1a^{n+1}}
\\
\cdot\cdot & \cdots & \cdot\cdot & \cdot\cdot
& \cdot\cdot & \cdots & \cdot\cdot
\\
Q_{x^ka^1} & \cdots & Q_{x^ka^{\mu-1}} & 1 & Q_{x^ka^{\mu+1}} 
& \cdots & Q_{x^ka^{n+1}}
\\
\cdot\cdot & \cdots & \cdot\cdot & \cdot\cdot
& \cdot\cdot & \cdots & \cdot\cdot
\\
Q_{x^na^1} & \cdots & Q_{x^na^{\mu-1}} & 0 & Q_{x^na^{\mu+1}} 
& \cdots & Q_{x^na^{n+1}}
\end{array}
\right\vert.
\endaligned
\]
To avoid any ambiguity, we shall sometimes put the integer $\mu$ in
the upper index position of the vertical bar to indicate precisely which
column is concerned. As is clear, this notation allows one to view and
to remember what are the involved partial derivatives of the
fundamental function $Q$ that appear inside each column. In summary,
$\square_{ [0_{ 1+ \ell}]}^\mu$ comes from $\square$ by changing just
its $\mu$-th column, as Cramer's rule classically says.

Next, the two-ways transfer between local functions $G$ defined in the
$(x, y, y_x)$-space and local functions $T$ defined in the $(x, a,
b)$-space, namely the one-to-one correspondence:
\[
G\big(x^1,\dots,x^n,y,y_{x^1},\dots,y_{x^n}\big)
\longleftrightarrow
T\big(x^1,\dots,x^n,a^1,\dots,a^n,b\big)
\] 
through the diffeomorphism~\thetag{ \ref{x-y-a-b}}, may be viewed
concretely, in the direction we are interested in, as the following
identity:
\[
\aligned
&
G\big(x^1,\dots,x^n,y,y_{x^1},\dots,y_{x^n}\big)
\equiv
\\
\equiv
T
\big(
&
x^1,\dots,x^n,\,
\\
&
A^1\big(x^1,\dots,x^n,y,y_{x^1},\dots,y_{x^n}\big),
\dots,
A^n\big(x^1,\dots,x^n,y,y_{x^1},\dots,y_{x^n}\big),
\\
&
A^{n+1}\big(x^1,\dots,x^n,y,y_{x^1},\dots,y_{x^n}\big)\big)
\endaligned
\]
holding of course in $\C \big\{ x^1, \dots, x^n, y, y_{ x^1}, \dots,
y_{ x^n} \big\}$. We therefore readily deduce how the derivation
$\frac{ \partial}{ \partial y_{ x^\ell}}$ is transferred to the $(x,
a, b)$-space:
\[
\frac{\partial G}{\partial y_{x^\ell}}
=
{\textstyle{\frac{\partial A^1}{\partial y_{x^\ell}}}}\cdot
\frac{\partial T}{\partial a^1}
+\cdots+
{\textstyle{\frac{\partial A^n}{\partial y_{x^\ell}}}}\cdot
\frac{\partial T}{\partial a^n}
+
{\textstyle{\frac{\partial A^{n+1}}{\partial y_{x^\ell}}}}\cdot
\frac{\partial T}{\partial a^{n+1}}. 
\]
By applying twice any two such derivations $\partial \big/ \partial
y_{ x^{ \ell_1 }}$ and $\partial \big/ \partial y_{ x^{ \ell_2} }$ to
an arbitrary function $G$, we may see,
after a few computations, what such a
composed differentiation corresponds to, in terms of the function $T$
defined in the $(x, a, b)$-space:
\[
\aligned
\frac{\partial^2G}{\partial
y_{x^{\ell_1}}\partial y_{x^{\ell_2}}}
&
=
\bigg(
\sum_{\mu=1}^{n+1}\,
{\textstyle{\frac{\partial A^\mu}{\partial y_{x^{\ell_1}}}}}\,
\frac{\partial}{\partial a^\mu}
\bigg)
\bigg[
\sum_{\nu=1}^{n+1}\,
{\textstyle{\frac{\partial A^\nu}{\partial y_{x^{\ell_2}}}}}\,
\frac{\partial T}{\partial a^\nu}
\bigg]
\\
&
=
\sum_{\mu=1}^{n+1}\,\sum_{\nu=1}^{n+1}\,
{\textstyle{\frac{\partial A^\mu}{\partial y_{x^{\ell_1}}}}}\,
{\textstyle{\frac{\partial A^\nu}{\partial y_{x^{\ell_2}}}}}\,
\frac{\partial^2T}{\partial a^\mu\partial a^\nu}
+
\sum_{\mu=1}^{n+1}\,\sum_{\nu=1}^{n+1}\,
{\textstyle{\frac{\partial A^\mu}{\partial y_{x^{\ell_1}}}}}\,
{\textstyle{\frac{\partial}{\partial a^\mu}}}
\big[
{\textstyle{\frac{\partial A^\nu}{\partial y_{x^{\ell_2}}}}}
\big]
\frac{\partial T}{\partial a^\nu}.
\endaligned
\]
Here, by a helpful formal convention, the three Greek letters $\mu$,
$\nu$ and $\tau$ will be used as summation indices in the total set
$\{ 1, \dots, n, n+1 \}$, while the four Latin letters $i$, $j$, $k$,
$\ell$ will always run in the restricted set $\{ 1, \dots, n\}$.
Replacing then the partial derivatives of the $A^\mu$ by their
values~\thetag{ \ref{A-square} } obtained previously, we thus get:
\[
\footnotesize
\aligned
\frac{\partial^2G}{\partial
y_{x^{\ell_1}}\partial y_{x^{\ell_2}}}
&
=
\sum_{\mu=1}^{n+1}\,\sum_{\nu=1}^{n+1}\,
\frac{\square_{[0_{1+\ell_1}]}^\mu}{\square}\,
\frac{\square_{[0_{1+\ell_2}]}^\nu}{\square}\,
\frac{\partial^2T}{\partial a^\mu\partial a^\nu}
+
\\
&
\ \ \ \ \ 
+
\sum_{\mu=1}^{n+1}\,\sum_{\nu=1}^{n+1}\,
\left\{
\frac{\square_{[0_{1+\ell_1}]}^\mu}{\square}\cdot
\frac{
\square\cdot
\frac{\partial}{\partial a^\mu}
\big(\square_{[0_{1+\ell_2}]}^\nu\big)
-
\square_{[0_{1+\ell_2}]}^\nu\cdot
\frac{\partial}{\partial a^\mu}
\big(\square\big)
}{\square\cdot\square}
\right\}
\frac{\partial T}{\partial a^\nu}
\endaligned
\]
Here, the coefficients of the $\frac{ \partial^2 T}{ \partial a^\mu
\partial a^\nu}$ will not be touched anymore, but the coefficients of
the $\frac{ \partial T}{ \partial a^\nu}$ must be subjected to further
transformations towards formal harmony, especially the numerator
involving a subtraction.

First of all, let us rewrite in length the concerned partial
derivative of the appearing modified Jacobian determinant\footnote{\,
Remind that, in order to differentiate a determinant, one should
differentiate separately each column once and then sum all the
obtained terms.  
}: 
\[
\aligned
{\textstyle{\frac{\partial}{\partial a^\mu}}}
\big(\square_{[0_{1+\ell_2}]}^\nu\big)
&
=
{\textstyle{\frac{\partial}{\partial a^\mu}}}
\Big[
\square
\big(
a^1\vert\cdots\vert a^{\nu-1}\vert 0_{[1+\ell_2]}\vert a^{\nu+1}
\vert\cdots\vert a^{n+1}
\big)
\Big]
\\
&
=
\square
\big(
a^1a^\mu\vert\cdots\vert a^{\nu-1}\vert 0_{[1+\ell_2]}\vert 
a^{\nu+1}\vert\cdots\vert a^{n+1}
\big)
+\cdots+
\\
&
\ \ \ \ \
+
\square
\big(
a^1\vert\cdots\vert a^{\nu-1}a^\mu\vert 0_{[1+\ell_2]}\vert 
a^{\nu+1}\vert\cdots\vert a^{n+1}
\big)
\\
&
\ \ \ \ \
+0+
\\
&
\ \ \ \ \
+
\square
\big(
a^1\vert\cdots\vert a^{\nu-1}\vert 0_{[1+\ell_2]}\vert 
a^{\nu+1}a^\mu\vert\cdots\vert a^{n+1}
\big)
+\cdots+
\\
&
\ \ \ \ \
+
\square
\big(
a^1\vert\cdots\vert a^{\nu-1}\vert 0_{[1+\ell_2]}\vert 
a^{\nu+1}\vert\cdots\vert a^{n+1}a^\mu
\big),
\endaligned
\] 
and also at the same time the partial derivative of the plain
Jacobiant determinant:
\[
\aligned
{\textstyle{\frac{\partial}{\partial a^\mu}}}
\big(\square\big)
&
=
{\textstyle{\frac{\partial}{\partial a^\mu}}}
\Big[
\square
\big(
a^1\vert\cdots\vert a^{\nu-1}\vert a^\nu\vert a^{\nu+1}\vert
\cdots\vert a^{n+1}
\big)
\Big]
\\
&
=
\square
\big(
a^1a^\mu\vert\cdots\vert a^{\nu-1}\vert a^\nu\vert a^{\nu+1}
\vert\cdots\vert a^{n+1}
\big)
+\cdots+
\\
&
\ \ \ \ \
+
\square
\big(
a^1\vert\cdots\vert a^{\nu-1}a^\mu\vert a^\nu\vert a^{\nu+1}
\vert\cdots\vert a^{n+1}
\big)
+
\\
&
\ \ \ \ \
+
\square
\big(
a^1\vert\cdots\vert a^{\nu-1}\vert a^\nu a^\mu\vert a^{\nu+1}
\vert\cdots\vert a^{n+1}
\big)
+
\\
&
\ \ \ \ \
+
\square
\big(
a^1\vert\cdots\vert a^{\nu-1}\vert a^\nu\vert a^{\nu+1}a^\mu
\vert\cdots\vert a^{n+1}
\big)
+\cdots+
\\
&
\ \ \ \ \
+
\square
\big(
a^1\vert\cdots\vert a^{\nu-1}\vert a^\nu\vert a^{\nu+1}
\vert\cdots\vert a^{n+1}a^\mu
\big). 
\endaligned
\]
Consequently, the numerator with a subtraction that we want to
simplify may be rewritten in length as follows:
\begin{equation}
\label{numerator-subtract}
\footnotesize
\aligned
\square\cdot
{\textstyle{\frac{\partial}{\partial a^\mu}}}
\big(\square_{[0_{1+\ell_2}]}^\nu\big)
-
\square_{[0_{1+\ell_2}]}^\nu\cdot
{\textstyle{\frac{\partial}{\partial a^\mu}}}
\big(\square\big)
=
\endaligned
\end{equation}
\[
\footnotesize
\aligned
&
=
\underline{\square
\big(
a^1\vert\cdots\vert a^{\nu-1}\vert a^\nu\vert a^{\nu+1}
\vert\cdots\vert a^{n+1}
\big)
\cdot
\square
\big(
a^1a^\mu\vert\cdots\vert a^{\nu-1}\vert 0_{[1+\ell_2]}\vert
a^{\nu+1}\vert\cdots\vert a^{n+1}
\big)}_{
\tiny{\octagon\!\!\!\! \sf a}}
+\cdots+
\\
&
\ \ \ \ \
+
\underline{\square
\big(
a^1\vert\cdots\vert a^{\nu-1}\vert a^\nu\vert a^{\nu+1}
\vert\cdots\vert a^{n+1}
\big)
\cdot
\square
\big(
a^1\vert\cdots\vert a^{\nu-1}a^\mu\vert 0_{[1+\ell_2]}\vert
a^{\nu+1}\vert\cdots\vert a^{n+1}
\big)}_{
\tiny{\octagon\!\!\!\! \sf b}}
+
\\
&
\ \ \ \ \
+0+
\\
&
\ \ \ \ \
+
\underline{\square
\big(
a^1\vert\cdots\vert a^{\nu-1}\vert a^\nu\vert a^{\nu+1}
\vert\cdots\vert a^{n+1}
\big)
\cdot
\square
\big(
a^1\vert\cdots\vert a^{\nu-1}\vert 0_{[1+\ell_2]}\vert
a^{\nu+1}a^\mu\vert\cdots\vert a^{n+1}
\big)}_{
\tiny{\octagon\!\!\!\! \sf c}}
+\cdots+
\\
&
\ \ \ \ \
+
\underline{\square
\big(
a^1\vert\cdots\vert a^{\nu-1}\vert a^\nu\vert a^{\nu+1}
\vert\cdots\vert a^{n+1}
\big)
\cdot
\square
\big(
a^1\vert\cdots\vert a^{\nu-1}\vert 0_{[1+\ell_2]}\vert
a^{\nu+1}\vert\cdots\vert a^{n+1}a^\mu
\big)}_{
\tiny{\octagon\!\!\!\! \sf d}}
-
\endaligned
\]
\[
\footnotesize
\aligned
&
-
\underline{\square
\big(
a^1\vert\cdots\vert a^{\nu-1}\vert 0_{[1+\ell_2]}\vert
a^{\nu+1}\vert\cdots\vert a^{n+1}
\big)
\cdot
\square
\big(
a^1a^\mu\vert\cdots\vert a^{\nu-1}\vert a^\nu\vert a^{\nu+1}
\vert\cdots\vert a^{n+1}
\big)}_{
\tiny{\octagon\!\!\!\! \sf a}}
-\cdots-
\\
&
-
\underline{\square
\big(
a^1\vert\cdots\vert a^{\nu-1}\vert 0_{[1+\ell_2]}\vert
a^{\nu+1}\vert\cdots\vert a^{n+1}
\big)
\cdot
\square
\big(
a^1\vert\cdots\vert a^{\nu-1}a^\mu\vert a^\nu\vert a^{\nu+1}
\vert\cdots\vert a^{n+1}
\big)}_{
\tiny{\octagon\!\!\!\! \sf b}}
-
\\
&
-
\underline{\square
\big(
a^1\vert\cdots\vert a^{\nu-1}\vert 0_{[1+\ell_2]}\vert
a^{\nu+1}\vert\cdots\vert a^{n+1}
\big)
\cdot
\square
\big(
a^1\vert\cdots\vert a^{\nu-1}\vert a^\nu a^\mu\vert a^{\nu+1}
\vert\cdots\vert a^{n+1}
\big)}_{\sf OK}
-
\\
&
-
\underline{\square
\big(
a^1\vert\cdots\vert a^{\nu-1}\vert 0_{[1+\ell_2]}\vert
a^{\nu+1}\vert\cdots\vert a^{n+1}
\big)
\cdot
\square
\big(
a^1\vert\cdots\vert a^{\nu-1}\vert a^\nu \vert a^{\nu+1}a^\mu
\vert\cdots\vert a^{n+1}
\big)}_{
\tiny{\octagon\!\!\!\! \sf c}}
-\cdots-
\\
&
-
\underline{\square
\big(
a^1\vert\cdots\vert a^{\nu-1}\vert 0_{[1+\ell_2]}\vert
a^{\nu+1}\vert\cdots\vert a^{n+1}
\big)
\cdot
\square
\big(
a^1\vert\cdots\vert a^{\nu-1}\vert a^\nu \vert a^{\nu+1}
\vert\cdots\vert a^{n+1}a^\mu
\big)}_{
\tiny{\octagon\!\!\!\! \sf d}}.
\endaligned
\]
The ante-penultimate underlined term ``${\scriptstyle\sf OK}$'' will
be kept untouched. To the pairs of (subtracted) $\square$-binomials
that are underlined with ${\sf a}$, ${\sf b}$, ${\sf c}$, ${\sf d}$
appended (including of course all terms present in the four
``$\cdots$''), we need an elementary instance of the Plücker
identities.

To state it generally, let $m \geqslant 2$, let $C_1, C_2, \dots, C_m,
D, E$ be $(m+2)$ column vectors in $\K^m$ and introduce the following
notation for the $m\times (m+2)$ matrix consisting of these vectors:
\[
\left[
C_1 \vert C_2 \vert \cdots \vert
C_m \vert D \vert E 
\right].
\]
Extracting columns from this matrix, we shall construct $m\times m$
determinants that are modification of the following ``ground''
determinant:
\[
\left\vert \!
\left\vert
C_1 \vert \cdots \vert
C_m
\right\vert \!
\right\vert\equiv
\left\vert \!
\left\vert
C_1 \vert \cdots \vert^{j_1}C_{j_1} \vert
\cdots \vert^{j_2} C_{j_2} \vert \cdots \vert
C_m
\right\vert \!
\right\vert.
\]
We use a double vertical line in the beginning
and in the end to denote a determinant. Also, we emphasize two
distinct columns, the $j_1$-th and the $j_2$-th, where $j_2 > j_1$,
since we will modify them. For instance in this matrix, let us replace
these two columns by the column $D$ and by the column $E$, which
yields the determinant:
\[
\left\vert \!
\left\vert
C_1 \vert \cdots \vert^{j_1} D \vert
\cdots \vert^{j_2} E \vert \cdots \vert
C_m
\right\vert \!
\right\vert.
\]
In this notation, one should understand that {\it only}\, the $j_1$-th
and the $j_2$-th columns are distinct from the columns of the
fundamental $m\times m$ ``ground'' determinant.

\smallskip\noindent{\bf Lemma.}
(\cite{ me2006}, p.~155)
{\em 
The following quadratic identity between determinants holds
true:}
\[
\aligned
& \
\left\vert \! 
\left\vert
C_1 \vert \cdots \vert^{j_1} D 
\vert \cdots \vert^{j_2}
E \vert \cdots \vert C_n
\right\vert \! 
\right\vert
\cdot
\left\vert \! 
\left\vert
C_1 \vert \cdots \vert^{j_1} C_{j_1}
\vert \cdots \vert^{j_2}
C_{j_2} \vert \cdots \vert C_n
\right\vert \! 
\right\vert = \\
& \
=
\left\vert \! 
\left\vert
C_1 \vert \cdots \vert^{j_1} D 
\vert \cdots \vert^{j_2}
C_{j_2} \vert \cdots \vert C_n
\right\vert \! 
\right\vert
\cdot
\left\vert \! 
\left\vert
C_1 \vert \cdots \vert^{j_1} C_{j_1}
\vert \cdots \vert^{j_2}
E \vert \cdots \vert C_n
\right\vert \! 
\right\vert- \\
& \
-
\left\vert \! 
\left\vert
C_1 \vert \cdots \vert^{j_1} E 
\vert \cdots \vert^{j_2}
C^{j_2} \vert \cdots \vert C_n
\right\vert \! 
\right\vert
\cdot
\left\vert \! 
\left\vert
C_1 \vert \cdots \vert^{j_1} C_{j_1}
\vert \cdots \vert^{j_2}
D \vert \cdots \vert C_n
\right\vert \! 
\right\vert.
\endaligned
\]

Admitting this elementary statement without redoing its proof and
applying it to all the above underlined pairs of (subtracted)
monomials, after checking that all final signs are ``$-$'', we obtain
the following neat expression for~\thetag{ \ref{numerator-subtract}}:
\[
\aligned
&
\square\cdot
{\textstyle{\frac{\partial}{\partial a^\mu}}}
\big(\square_{[0_{1+\ell_2}]}^\nu\big)
-
\square_{[0_{1+\ell_2}]}^\nu\cdot
{\textstyle{\frac{\partial}{\partial a^\mu}}}
\big(\square\big)
=
\\
&
=
-
\square
\big(
0_{[1+\ell_2]}\vert\cdots\vert^\nu\,a^\nu\vert
\cdots\vert a^{n+1}
\big)
\cdot
\square
\big(
a^1\vert\cdots\vert^\nu\,a^1a^\mu\vert
\cdots\vert a^{n+1}
\big)
-\cdots-
\\
&
\ \ \ \ \
-
\square
\big(
a^1\vert\cdots\vert^\nu\,0_{[1+\ell_2]}\vert
\cdots\vert a^{n+1}
\big)
\cdot
\square
\big(
a^1\vert\cdots\vert^\nu\,a^\nu a^\mu\vert
\cdots\vert a^{n+1}
\big)
-\cdots-
\\
&
\ \ \ \ \
-
\square
\big(
a^1\vert\cdots\vert^\nu\,a^\nu\vert
\cdots\vert 0_{[1+\ell_2]}
\big)
\cdot
\square
\big(
a^1\vert\cdots\vert^\nu\,a^{n+1}a^\mu\vert
\cdots\vert a^{n+1}
\big), 
\endaligned
\]
or equivalently, in contracted form: 
\[
\square\cdot
{\textstyle{\frac{\partial}{\partial a^\mu}}}
\big(\square_{[0_{1+\ell_2}]}^\nu\big)
-
\square_{[0_{1+\ell_2}]}^\nu\cdot
{\textstyle{\frac{\partial}{\partial a^\mu}}}
\big(\square\big)
=
-\sum_{\tau=1}^{n+1}\,
\square_{[0_{1+\ell_2}]}^\tau
\cdot
\square_{[a^\tau a^\mu]}^\nu. 
\]
Thanks to this sidework, coming back to the expression for $\frac{
\partial^2 G}{ \partial y_{ x^{\ell_1}} \partial y_{ x^{\ell_2}}}$ we
left pending above, we obtain:
\[
\aligned
\frac{\partial^2G}{\partial y_{x^{\ell_1}}\partial y_{x^{\ell_2}}}
&
=
\frac{1}{\square^2}\,
\sum_{\mu=1}^{n+1}\,\sum_{\nu=1}^{n+1}\,
\left\{\,
\square_{[0_{1+\ell_1}]}^\mu
\!\cdot
\square_{[0_{1+\ell_2}]}^\nu
\right\}
\frac{\partial^2T}{\partial a^\mu\partial a^\nu}
-
\\
&
\ \ \ \ \ 
-
\frac{1}{\square^3}\,
\sum_{\mu=1}^{n+1}\,
\sum_{\nu=1}^{n+1}\,
\sum_{\tau=1}^{n+1}\,\,
\left\{
\square_{[0_{1+\ell_1}]}^\mu
\cdot
\square_{[0_{1+\ell_2}]}^\tau
\cdot
\square_{[a^\mu a^\tau]}^\nu
\right\}
\frac{\partial T}{\partial a^\nu}.
\endaligned
\]
To really finalize this expression, we factor everything by $\frac{
1}{ \square^3}$ and we exchange the two summation indices $\nu$ and
$\tau$ in the second line:
\[
\boxed{
\aligned
\frac{\partial^2G}{\partial y_{x^{\ell_1}}\partial y_{x^{\ell_2}}}
&
=
\frac{1}{\square^3}\,
\sum_{\mu=1}^{n+1}\,\sum_{\nu=1}^{n+1}\,
\square_{[0_{1+\ell_1}]}^\mu
\!\cdot
\square_{[0_{1+\ell_2}]}^\nu
\left\{\,
\square\cdot
\frac{\partial^2T}{\partial a^\mu\partial a^\nu}
-
\sum_{\tau=1}^{n+1}\,\,
\square_{[a^\mu a^\nu]}^\tau
\cdot
\frac{\partial T}{\partial a^\tau}
\right\}
\endaligned\!\!\!\!}\,.
\]

\smallskip\noindent
{\em End of proof of the Main Theorem.}  As already explained in the
Introduction, one applies to the system~\thetag{ \ref{Phi-second}}
Hachtroudi's characterization~\thetag{ \ref{zero-hachtroudi}} of
equivalence to the system $w_{ z_{ k_1}' z_{ k_2}'} ' (z') = 0$ with
$x := z$, with $y := w$, with $a := \overline{ z}$, with $b :=
\overline{ w}$, with $(a, b) := \overline{ t}$, with $Q := \Theta$,
with $\square := \Delta$, with $G := \Phi_{ k_1, k_2}$ and with $T :=
\frac{ \partial^2 \Theta}{ \partial z_{ k_1} \partial z_{ k_2}}$.  The
denominator $\frac{ 1}{ \Delta^3}$ can be cleared out, and we simply
get the explicit fourth-order partial differential equation satisfied
by $\Theta$. This completes the proof of our Main Theorem
and the paper may end up now.
\qed

\vfill\end{document}